\documentclass[11pt]{article}
\usepackage[utf8]{inputenc}
\usepackage[
papersize={5.5in, 8.5in},
left=0.35in,
right=0.35in,
top=0.35in,
bottom=0.5in]{geometry}
\usepackage{graphicx}
\usepackage{amsthm}
\usepackage{rotating}
\usepackage{amsfonts,amssymb,amsmath}

\usepackage{booktabs}
\usepackage[table,xcdraw]{xcolor}

\newcommand{\gr}{\cellcolor[HTML]{E0E0E0}}

\DeclareMathOperator{\Aut}{Aut}
\graphicspath{ {images/} }

\title{\bf 
Graph minimization,\\ focusing on the example
of 5-chromatic\\ unit-distance graphs in the plane}
\author{\bf \normalsize Jaan Parts} 
\date{\normalsize Kazan, Russia, jaan\_parts@mail.ru}

\begin{document}

\maketitle

\pagestyle{empty}
\thispagestyle{empty}

\begin{abstract}
We introduce a new graph minimization method, in which it is required to preserve some graph property and there is an effective procedure for checking this property. We applied this method to minimize 5-chromatic unit-distance graphs and obtained a graph with 509 vertices and 2442 edges.
\end{abstract}

\section{Hi story}

Everybody repeats this fabulous story. I want to tell it too. 

Once upon a time, an old wizard (good or evil, nobody knows) gave birth to a multi-headed mathematical dragon. In fact, all these heads, except for one, were superfluous
(the secret is that this is not a dragon, but an enchanted princess), but the question is which head is real? (He who answers will save the princess and get half the kingdom.) Since the dragon is mathematical, all his heads are numbered. Almost all of them fell off right away, but the four that survived (numbered 4, 5, 6, and 7) held tight and were extremely viable.

The witchcraft spell that created the dragon was a simple question: how many colors are needed to color a plane so that any two points at a unit distance from each other have a different color? (In magic language, this number is called the \textit{chromatic number} and is denoted by the magic sign~$\chi$.)

Many heroes went to battle with the dragon, but their spears broke, the swords were dulled, and the brave men were forced to retreat. While they thought how to deal with extra heads, and examined this problem from different angles, the dragon began to breed. A great number of small dragonlings appeared, and almost all of them inherited many heads (and the more distant the descendant, the more heads he has). From time to time it turned out to be possible to remove one or two, but the dragons continue to breed.

But what about the main dragon? For thirty years and another thirty-three years (or so) no one knew how to approach him. Regarding which head is real, various considerations and predictions have been expressed. The storyteller saw a mathematician who claimed that, well, almost everyone is sure that the real head is the one with the number 4. A week later it was proved that this head was superfluous. Aubrey came, and just chopped it off.

But 33.3\% of success\footnote{It is assumed that the goal is to chop off all but one head. Although it is possible that the wizard decided to play a trick and enchanted the Siamese twins.}
was not enough for him, so he assembled a team (called the Polymath project), and all together they rushed to pursue a fast-fading monster. But the latter managed to slip away. This is not the end of the story, and it continues today. Maybe the dragon is hiding behind the next bend, and one just has to look carefully, and the head number 5 disappears by itself?

So far, some battlefield-cleaning tasks remain. One of them is to clarify the size of the dragon, to study the moment of separation of head number~4. This can provide a clue to understanding how to proceed. And here we are already forced to completely set aside fabulous terminology and move closer to mathematics.

\section{Toolkit}

\paragraph{Unit distance graphs.}
To prove that $\chi\ge5$ is true for the plane, it is enough to construct a 5-chromatic unit distance graph.

Formally, a \textit{graph} is a pair $G =(V, E)$, where $V$ is a collection of elements called \textit{vertices}, $E$ is a collection of two-element subsets of $V$ called \textit{edges}. A special class of graphs is the \textit{distance graphs}, which can be \textit{embedded}, i.e. each vertex can be assigned coordinates in $n$-dimensional space for some specified $n$ (in our case $n=2$), in such a way that each edge connects only pairs of vertices that are at certain Euclidean distances from each other: in our case, at unit distance. We restrict ourselves in this paper to the consideration of \textit{strict unit distance graphs} for which all pairs of vertices located at unit distance are connected by an edge. (This restriction is not necessary for the approach we describe, but simplifies the presentation and comparison of the results.) In this case, the graph can be identified by the list of vertex coordinates, since the latter completely defines the sets $V$ and $E$.

\paragraph{Problem formulation.}
Our main focus is to minimize the number of vertices of a unit distance graph while maintaining the chromatic number 5. Since one can usually find many non-isomorphic graphs with a given number of vertices, it is customary to choose among these the graph with the fewest edges.

The main objective of this paper is to describe a new method for minimizing graphs while preserving a desired property (in particular, the chromatic number).

\paragraph{Definitions.}
A \textit{proper vertex coloring} is a way of coloring graph vertices such that no two adjacent vertices (in our case at a unit distance apart) are of the same color. If $k$ colors are sufficient, then the graph is called $k$-\textit{colorable} (below default $k=4$). The minimum number $k$ such that a given graph is $k$-colorable is the graph’s \textit{chromatic number} $\chi$.

A \textit{monochromatic pair} of vertices (or a \textit{mono-pair}) are two vertices of a particular graph that have the same color in any proper $k$-coloring. 
Similarly, a \textit{non-monochromatic pair} (or \textit{non-mono-pair}) consists of two vertices that cannot have the same color (a trivial example is a pair that defines an edge).

A non-mono pair that does not define an edge has sometimes \cite{pm16} been termed a \textit{virtual edge}, but in our opinion, this term is not entirely successful, since it does not clearly express whether one means a mono- or non-mono-pair, which leads to confusion. 

A cycle of two mono-pairs and a unit edge is called a \textit{spindle}.

A \textit{non-monochromatic triple} (or \textit{non-mono-triple}) is a set of three vertices that cannot simultaneously have the same color. If the triple of vertices forms an equilateral triangle, then it is characterized by the length of the \textit{side}.

A graph having a certain subgraph is, for the latter, an \textit{overgraph}.

An \textit{automorphism} of a graph $G$ is a permutation on the set of vertices that preserves the adjacency relation. The set of all permutations forms a group $\Aut G$ of automorphisms (or \textit{symmetries}). The set of vertices is divided into a disjoint set of \textit{orbits}, which are equivalence classes defined as follows: two vertices belong to the same orbit of the graph $G$ if and only if there is an automorphism group $\Aut G$ that maps one of them to the other.

\paragraph{Basic operations.}
At the moment, truly small 5-chromatic graphs are unknown: one needs to operate with hundreds and thousands of vertices. Therefore, an effective procedure is needed to create large structured graphs. The basic tools for this are (a) Minkowski summing and (b) rotations that allow the addition of new edges.

The \textit{Minkowski sum} $G=G_1\oplus G_2$ of two unit-distance graphs $G_1=(V_1, E_1)$ and $G_2=(V_2, E_2)$ is a graph $G=(V, E)$ whose vertex set is determined by the union of all possible sums of their vertex vectors $V=\{v\,|\,v=v_1+v_2,v_1\in V_1, v_2\in V_2\}$, and the set of edges includes all obtained vertex pairs with unit distance $E=\left\{\{v_a,v_b\}\,|\,|v_a-v_b|=1,\,v_a, v_b\in V \right\}$.

\paragraph{Calculation of graph chromatic number.}
The task of minimizing 5-chromatic unit-distance graphs would be difficult to accomplish if there were no effective computational procedure for determining the chromatic number. Fortunately, such a procedure exists and uses the \textit{satisfiability} approach \cite{knu} implemented in SAT solvers, such as {\tt MiniSAT} (also integrated into {\tt Mathematica}), {\tt Glucose} and others. An equivalent task is to check whether a given graph is 4-colorable or not.

\paragraph{Notation.}
For short, we denote the set with elements $a_n$ by $\{a\}$, omitting the index. 
Graphs are denoted by a capital italic letter, for example $G$. Most letters correspond to special graphs. 
The graph symbol may be preceded by a complex rotation and (or) scaling multiplier denoted by a lowercase italic letter.
In particular, the following rotation multipliers are used: $i=\exp\frac{i\pi}{2}=\sqrt{-1}$, $\eta=\exp\left (\frac{i}{2}\arccos\frac56\right)$, $\rho=\exp\left (i\arccos\frac78\right)$. 

The subscript in a graph symbol indicates the number of its vertices. The superscript $\{a\}$ in $G^{\{a\}}$ defines a set of rotations by multiples of $\eta$: $G^{\{a\}}=\bigcup_{\alpha\in \{a\}}\eta^{\alpha}G$. In the special case $\{a\}=\{-m\ldots+m\}$, 
a simplified notation is 
used: $G^m= G^{\{-m\ldots+m\}}$. Note that $G^0= G$. The Minkowski sum of $n$ identical graphs $G$ we denote as $\oplus^n G$. Note that $\oplus G=G$.

The main building block is the 7-vertex hexagonal \textit{wheel} graph $H=H^0=W_7$ with vertex coordinates $(0,0)$, $(\pm 1,0)$, and $(\pm 1/2,\pm \sqrt3/2)$.
The 31-vertex graph used in \cite{grey, heu553} to construct 5-chromatic graphs is denoted as $ V_{31}=H^{\{-2,-1,\,0,\,1,\,2\}}=H^2$. Many of the graphs discussed below are subgraphs of a 12937-vertex graph $H^2\oplus H^2\oplus H^1\oplus H^1=\oplus^2H^2\oplus^2H^1$.

\paragraph{Base graph.}

Of particular interest is a construction of the form $\oplus^nH^m$, which we call the \textit{base graph} $B$. Here $n$ is the radius\footnote{For base graph the radius takes the same value, half the diameter, in both the graph-theoretic and the Euclidean sense.} of the graph, while $m$ is the number of rotations by $\eta$. In the complex plane, 
vertices of $B$ have coordinates of the form \cite{exoo, gib}
$$\frac{a+b\sqrt{33}+ic\sqrt{3}+id\sqrt{11}}{4\cdot3^h},
\;\;a,b,c,d \in \mathbb{Z}, \;a-b+c+d \in 4\mathbb{Z},\; 
h=\lceil m/2 \rceil
$$

In the graphs discussed below, two rotations $m$ in either direction (i.e. a total of five copies) are sufficient; up to three rotations are used in \cite{exoo, heu553}.
Below, each vertex of the graph is denoted by the vector $(a, b, c, d)$, and it is assumed that $h=1$.

When we talk about orbits, we usually (but not always) mean the orbits of the base graph, or in short, \textit{base orbits}.
The vertices of the orbit are sorted by the number of nonzero and negative elements in notation, then in lexicographic order.
The orbit is denoted by the first vertex in the list of its vertices. All vertices of the orbit can be obtained from the first one using the identical transform $\tau_0$, positive $\tau_1$ and negative $\tau_2$ rotations by $2\pi/3$, reflections relative to the vertical $\tau_3$ and horizontal $\tau_4$ axes and conjugates\footnote{A conjugate is a binomial formed by negating the second term of a binomial. In our case, the conjugate of $x+y\sqrt{33}$ is $x-y\sqrt{33}$, where $x=a+i c\sqrt{3}$, $y=b+i d/\sqrt{3}$.}
$\tau_5$ by $\sqrt{33}$.

The corresponding coordinate transform matrices\footnote{Any vertex $V_\tau$ of the orbit can be obtained from some vertex $V$ of this orbit by successive application of several coordinate transformations, $V_\tau=\tau V$, $\tau=\tau_0[\tau_1|\tau_2][\tau_3][\tau_4][\tau_5]$, where the square brackets mean an optional factor. }
are of the form (empty cells represent zero values):

\vspace{1.5mm}
{
\renewcommand{\arraystretch}{0}
\begin{tabular}{@{}rrr@{}}
$\tau _0=\left[ \begin{smallmatrix}\;\:1\;&&&\\&\;\:1\;&&\\&&\;\:1\;&\\&&&\;\:1\; \end{smallmatrix} \right]$ &
$\tau _1=\frac12 \left[ \begin{smallmatrix}-1&&-3&\\&-1&&-1\\1&&-1&\\&3&&-1 \end{smallmatrix} \right]$ &
$\tau _2=\frac12 \left[ \begin{smallmatrix}-1&&3&\\&-1&&1\\-1&&-1&\\&-3&&-1 \end{smallmatrix} \right]$ \\
\rule{0pt}{5pt}\\
$\tau _3=\left[ \begin{smallmatrix}-1&&&\\&-1&&\\&&\;\:1\;&\\&&&\;\:1\; \end{smallmatrix} \right]$ &
$\tau _4=\;\;\;\left[ \begin{smallmatrix}\;\:1\;&&&\\&\;\:1\;&&\\&&-1&\\&&&-1 \end{smallmatrix} \right]$ &
$\tau _5=\;\;\; \left[ \begin{smallmatrix}\;\:1\;&&&\\&-1&&\\&&\;\:1\;&\\&&&-1 \end{smallmatrix} \right]$ \\
\end{tabular}

}
\vspace{1.5mm}

The symmetry group of the base graph of the form $\oplus^nH^m$ ($n\ge2$, $m\ge1$) has maximum order 24 and includes the dihedral group and the additional symmetry formed by the conjugation transform $\tau_5$. A fully filled orbit can be of order 1, 6, 12, or 24. The list of orbits is sorted by the order of its base orbit, then by the smaller of the two possible Euclidean radii.

\section{Evolution}
We skip the entire previous evolution of coloring, which is well described in amazing Soifer's book \cite{soi}, and immediately move on to a short history of minimizing 5-chromatic unit-distance graphs in the plane.

The initial work on 5-chromatic graph minimization was done by Aubrey de Grey \cite{grey} himself soon after finding the first such graph with 20425 vertices. He obtained a construction of four 397-vertex identical subgraphs, which could be connected to give a 1585-vertex graph.

The main work on further minimization was done by Marijn Heule \cite{heu553}, who found a more effective way of joining subgraphs, reducing their number to two, moreover one of which can be much smaller than the other. Hereafter such graphs will be referred to as the combination of a large $L$ and a small $S$ subgraph, with a specified relative rotation about a shared vertex, i.e. $G=L \cup \rho S$. Heule’s method of minimization is based on the use of a special {\tt DRAT-trim} program that extracts additional information (the “unsatisfiable core”) from a SAT solution. He quickly discovered a graph \cite{heu553} with 553 vertices (with 420- and 134-vertex subgraphs and one common vertex); a year later a 529-vertex graph \cite{heu529} was published (with 393- and 137-vertex subgraphs).

We approached the minimization problem from the other end, namely, as the problem of simplifying the structure of a 5-chromatic graph (or decomposing it into a minimum number of simple components). In particular, we used a 2-fold Minkowski sum in subgraphs (previously, a 3-fold sum was required, of the form $\oplus^3 H^2$ or $\oplus^2 H^2\oplus H^0$). For example, as one starting-point, we used graphs $L_{727}=V_{49}\oplus V_{37}$, $S_{361}=V_{31}\oplus V_{25}$, where $ V_{49}=H^2\cup\frac{i}{\sqrt3} H^1$, $ V_{37}= H^1\cup\frac{i}{\sqrt3} H^1$, $ V_{31}= H^1\cup{i}{\sqrt3} H^{\left\{-1,\,1\right\}}$, $ V_{25}= H^1\cup\frac{i}{\sqrt3} H^0$.

Then the obtained subgraphs were reduced first with a maximum symmetry order 24 to $L_{451}$ and $S_{181}$, then with a symmetry order 6 to $L_{418}$ and $S_{151}$. This gave a 568-vertex graph (it's not bad, but still far from a record). A further decrease in the symmetry order and the introduction of the stage of adding new vertices made it possible to compete with results of Heule.

Together, in the course of a short competition, we brought the number of vertices down to 525, 517 (Heule) and, finally, 510 (Heule soon corrected his method and came to a similar result \cite{heu510}). Later, when preparing this text, we obtained a graph with 509 vertices.

\section{Graph types}

We can distinguish a number of currently known families of 5-chromatic graphs, differing in the method of construction. We call such a family a \textit{type}. The main element of construction is a base graph, which can be entirely placed in $\mathbb{Q}[\sqrt3, \sqrt{11}]$ by some embedding.
Then, the type and subtypes of graphs in our classification can be defined by the number and relative position of the base graphs. The type name is given by the first letter of the name of the participant in the Polymath project, who was the first to pay attention to it by giving a concrete example.

\paragraph{Type A} (Aubrey de Grey \cite{grey})
is formed by at least four base graphs, each pair of which is connected in such a way as to form a mono-pair in any 4-coloring. Two or more mono-pairs and one non-mono-pair (in particular, unit edge) are connected in a cycle. 

This is historically the first type of 5-chromatic graphs. In the original, mono-pairs have a length of 4 and are obtained by rotating one of the base graphs relative to the other by $\rho$.

Type A has the most transparent construction in terms of proving the lower bound of the chromatic number of the plane, $\chi\ge5$. It relies on the existence of a non-mono triple forming an equilateral triangle with side $\sqrt3$. In \cite{grey} the graph with such a triple has 397 vertices; in \cite{heu278} this was reduced to 278 vertices. The demonstration that a graph of this order has a non-mono triple requires computer verification, and, apparently, is the smallest known obstacle to a complete (human-verified) proof that $\chi\ge5$. A similar triple with side $1/\sqrt3$ was used in an alternative proof described in more detail below \cite{exoo}.

\paragraph{Type M} (Marijn Heule \cite{heu553})
is formed by two base graphs rotated relative to each other so as to obtain many additional unit edges of two or more kinds: that is, additional edges which, together with the center of rotation, produce differently-shaped unequal triangles.

In the original \cite{heu553}, the rotation 
multiplier is $\rho$, which in addition to isosceles triangles with side 2 gives triangles with sides $\frac{\sqrt{11}}{2}\pm\frac{\sqrt3}{6}$. Unlike graphs of type A, the existence of triangles of the second kind is essential. Working constructions of type M with other rotations are not known yet.

Type M is of the greatest interest, because so far it gives 5-chromatic graphs with the least number of vertices (a little bit over 500). Graphs of type M with a small number of vertices are described by the expression $G=L\cup \rho S$ and are formed by the union of two subgraphs $L$ and $S$, which differ significantly (about three-fold) in the number of vertices.

\paragraph{Type J} (Jannis Harder \cite{har})
is formed by two or more base graphs, each of which contains a mono-pair in any 4-coloring. These base graphs are then connected so that two or more mono-pairs and one unit edge form a cycle, which gives a 5-chromatic graph.

In the original \cite{har}, the mono-pair has a length of $8/3$, and the corresponding graph has 745 vertices. Using a suitable base graph, one can obtain mono-pairs of different lengths, in particular, lengths of the form $n(8/3)^{m/2}$, where $n,m\in\mathbb{Z}$. To show this, it is enough to notice that a mono-pair of length $r/\sqrt3$ can be formed from two mono-pairs of some length $r$ by rotating one of them around a common vertex by $\eta$, which leaves both mono-pairs in $\mathbb{Q}[\sqrt3, \sqrt{11}]$.

Type J has the simplest construction, since it requires the least number of connections between base graphs.

An alternative proof for the lower bound of the chromatic number of a plane, proposed by Exoo and Ismailescu \cite{exoo}, is reduced to a mono-pair of length $8/3$, and contains intermediate constructions: non-mono-pairs of length $\sqrt{11/3}$ and non-mono-triples forming an equilateral triangle with side $1/\sqrt3$. For the latter ones, the corresponding graph has 627 vertices.

\paragraph{Type T} (Tam\'as Hubai \cite{hub})
is formed by connecting multiple base graphs into a common lattice.
Compared to type M, graphs of type T have many vertices common to different base graphs.

In the original \cite{hub}, a rhombic lattice with rotations by $\rho$ is formed. The corresponding graph has 6937 vertices and is defined (slightly differently than in the original) by $G=\oplus^3V_{37}$, where $V_{37}=H^1\cup \rho H^1$.

Type T allows the use of the simplest base subgraphs (if we count the number of rotations of the hexagonal lattice), but the required number of subgraphs is large and leads to large overall graphs (we managed to reduce the number of vertices only to about 5000).

\paragraph{Subtypes of type M.}

Graphs of type M substantially use additional connection edges between subgraphs $L$ and $S$. Additional edges either connect vertices that are remote from the center of rotation at a Euclidean distance $2$ (we call these reference edges) or ones at distances $\frac{\sqrt{11}}{2}\pm\frac{\sqrt3}{6}$ (we call these auxiliary edges). It is natural to ask, what is the minimum required number of these edges? 

\begin{figure}[!t]
{
\small
\centering
\begin{tabular}{@{}ccc@{}}
    \includegraphics[scale=0.33]{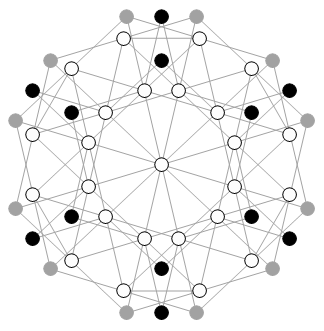} &
    \includegraphics[scale=0.33]{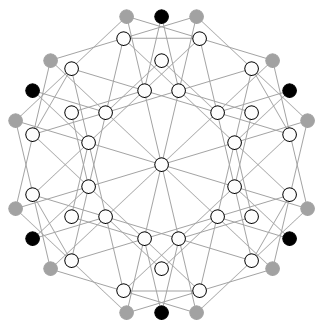} &
    \includegraphics[scale=0.33]{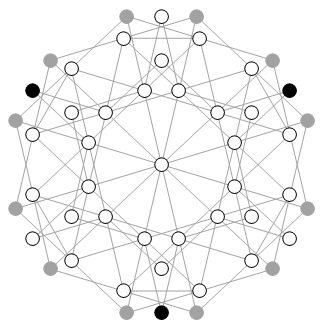} \\
    M12 & M6H & M3H \\
    \\
    \includegraphics[scale=0.33]{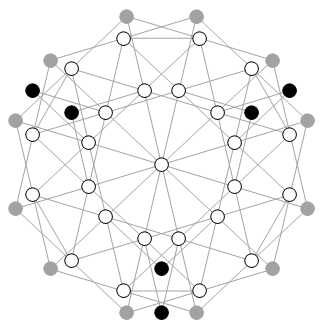} &
    \includegraphics[scale=0.33]{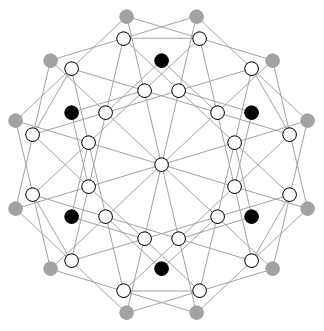} &
    \includegraphics[scale=0.33]{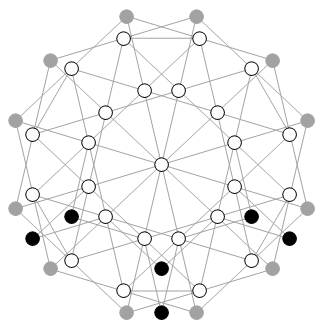} \\
    M6A & M6B & M6C \\
\end{tabular} \par
}
  \caption{
Some subtypes of type M. The vertices of one of the two subgraphs are shown. The reference orbit $(0,4,4,0)$ is shown in gray, the auxiliary orbit $(0,0,2,6)$ is shown in black. For clarity the orbits $(0,2,2,0)$ and $(6,0,0,6)$ are also shown.}
  \label{mtype}
\end{figure}

It was found that it is enough to leave two base orbits to connect $L$ and $S$, namely $(0,4,4,0)$ and $(0,0,2,6)$: see Fig.\ref{mtype}. Each of the orbits includes 12 vertices. The reference orbit $(0,4,4,0)$ is always full, while the auxiliary orbit $(0,0,2,6)$ may be only partially filled. The subgraphs $L$ and $S$ are connected by auxiliary edges crosswise. The number of auxiliary edges is included in the notation of a subtype of type M.

For relatively large graphs of type M, for example, of the form $L\cup \rho L$, three auxiliary edges are sufficient (subtype M3). If a small S-subgraph is used, then at least six auxiliary edges are necessary, and only three different (up to symmetry transformations) configurations are allowed. They are denoted M6A, M6B and M6C, and are shown at the bottom of Fig.\ref{mtype}. The subtype M6H is isomorphic to the subtype M6B. It can be seen that graphs of type M with larger numbers of auxiliary edges must include a subgraph of one of the three allowed subtypes M6.

The filled auxiliary orbits correspond to the subtype M12, the 553-vertex graph of Heule \cite{heu553} is of subtype M10, and the 529-vertex graph of Heule \cite{heu529} is of subtype M6A. The question arises, which subtype leads to the smallest 5-chromatic graph?

\section{Method}

\subsection{Preliminary considerations}
In general terms, the graph-minimization problem can be formulated as follows: having some initial graph with a specified “key property” (in our case, this property is the chromatic number of the graph or an overgraph of it), we should get a reduced graph that preserves this property and has the smallest possible number of vertices. Ideally, it is desirable to verify that a global minimum is obtained, at least within a limited class of graphs.

It makes sense to explicitly separate the concepts of reduction and minimization. By \textit{reduction}, we mean obtaining a subgraph of some given graph. By \textit{minimization}, we mean obtaining a graph with fewer vertices, but allowing the possibility of adding new vertices that were not in the initial graph.

Suppose that all that we have initially is a procedure for checking whether a graph possesses the key property, and on that basis we want to build a minimization procedure. Let's see what difficulties may arise.

The search for a minimal graph can be started by choosing a large graph that has the required property. The greater the order of the initial graph, the greater the likelihood that the minimal graph we seek is a subgraph of it. However, the complexity of the minimization problem increases. A global minimum can be achieved by checking all the options (brute force method), but usually the number of such options is huge. For example: suppose we reduce a graph containing 169 vertices and obtain a graph with 136 vertices. In a completely crazy version of the brute force method, we would check all ${169\choose 136}\approx 10^{35}$ subgraphs of order 136.

The opposite approach is to test each possible pair of vertices of the initial minimal graph by replacing them with one new vertex from some list of suitable ones. However, it was found that such a simple procedure, as a rule, does not allow us to get out of a local minimum. For this approach to work, one must operate with extremely large sets of vertices, which makes it computationally inefficient.

Therefore, the minimization procedure, on the one hand, should operate with a relatively small initial graph and be able to effectively reduce the number of checks, and on the other hand, it should not depend heavily on the choice of initial graph and should include the addition of new vertices to escape local minima.

\subsection{Main idea}
Our approach is based on a simple idea: if deleting a certain vertex (or set of vertices) of a graph leads to the loss of the key property, then there is no need to check any subgraphs that do not contain this vertex. This idea can be applied repeatedly, trying to remove one vertex first, then pairs, triples and so on until the time for checking the remaining subgraphs becomes comparable with the time spent on their further sifting.

\subsection{Variables}

Next, we will consider details of this approach, including the addition of new vertices, the use of graph symmetry and other ways to speed up calculations. But first, we list the graphs used in the minimization procedure:
\begin{itemize}
\itemsep=0pt
    \item
$M$ is the \textit{minimal} graph, which has the key property and has the smallest number of vertices found
    \item
$\{M\}$ is the \textit{set} of minimal graphs, i.e. all found graphs with the least number of vertices
    \item
$A$ is the \textit{accumulative} graph, which is the union of all the graphs of the set $\{M\}$
    \item
$W$ is the \textit{working} graph, which should be reduced at the next iteration of the minimization procedure and includes graph $A$ and some set of candidate new vertices 
    \item
$R$ is the \textit{reserve} graph, some high-symmetry-order overgraph of graph $A$, from which new vertices are selected to obtain the working graph $W$
    \item
$B$ is the \textit{base} graph of the search area, i.e. some overgraph of graph $R$, which includes all vertices to be considered during the minimization process
\end{itemize}

The listed graphs are related as follows:
$M\subseteq A\subseteq W\subseteq R\subseteq B$,
$M\in \{M\}$,
$\bigcup\{M\}=A$.

The following two graphs, generally speaking, are not unit distance graphs, and are introduced in order to describe used entities from a single perspective:

\begin{itemize}
\itemsep=0pt
    \item 
$C$ is the \textit{companion} graph, which is added to graph $W$ to obtain a property that can be conveniently checked using the SAT solver (for an $L$-graph, it is an $S$-graph and vice versa, for a mono-pair it is an edge of the corresponding length, for a non-mono-set it is 3-clique, each vertex of which is connected to all vertices of this set)
    \item 
$Y$ is a \textit{relationships} hypergraph, which contains vertex sets that cannot be wholly removed during reduction of graph $W$. 
A \textit{hypergraph} is a generalization of a graph in which an edge (called a \textit{hyperedge}) can join any number (called the \textit{hyperedge degree}) of vertices.
\end{itemize}

\subsection{Symmetries}

To speed up the calculations, we use graph symmetries, which in particular allow us to operate with vertex sets, or orbits, instead of individual vertices. In this case, it is necessary to distinguish:

\begin{itemize}
\itemsep=0pt
    \item 
base symmetry group $\Aut B$, which appears in base orbit filling tables (see section \ref{results})
    \item
the symmetry $\Aut W$ of working graph $W$ itself
    \item
the “input symmetry” $\Aut (W\cup C)$ 
at the start of the reduction procedure taking into account the companion graph $C$
    \item
the “output symmetry” $\Aut\{M\}$ of graph set $\{M\}$ at the end of the reduction procedure, i.e. the symmetry of smallest order $|\Aut M|:\forall M\in\{M\}$
    \item
the actual symmetry $\Aut M$ of a particular graph $M$, which may differ from $\Aut\{M\}$
\end{itemize}

\subsection{Algorithm}
In general, the minimization procedure involves the alternation of the operations of expansion and reduction of accumulative graph $A$, comprising the three following main stages.

1. In the expansion stage, several new vertices from the base graph $B$ 
are added to a graph $A$; as a result, a working graph $W$ is obtained, which then undergoes reduction.

2. In the first phase of the reduction, a hypergraph $Y$ is constructed, the hyperedges of which determine the critical relationships between the vertices of $W$. A hyperedge of degree $n$ is a subset of $n$ vertices of $W$ such that deleting all $n$ vertices at once leads to loss of the key property (in our case it means a decrease in the chromatic number of the graph $W \cup C$). 
Moreover, any proper subset of the vertices of a hyperedge should not lead to the loss of the key property when removed from $W$. In other words a hyperedge of some degree does not contain entirely hyperedges of lower degrees. 
Hyperedges are sought stepwise, in increasing order of their degree, up to a maximum that is selected so as to minimize the total calculation time.

3. In the second phase of reduction, a set of minimal graphs $\{M\}$ is found. To do this, the key property is checked for all possible subgraphs of the graph $W$ containing at least one vertex from each hyperedge of  $Y$. The verification is performed stepwise in order of increasing number of vertices of the subgraphs. At each step, all subgraphs of some order are checked. The check is terminated if at some step at least one subgraph with the key property is found. In any case, it is enough to check only the subgraphs whose order does not exceed the order of the minimum graph $M$ found at the previous iteration. The check can also be terminated if the number of checked subgraphs at some step exceeds a threshold.

Upon completion of the iteration of the minimization procedure, a new accumulative graph $A$ is formed, combining all the found minimal graphs $\{M\}$. The iteration is considered successful if it leads to a change in the graph $A$: that is, either to a reduction of the order of the minimum graph $M$ or else to an expansion of the set of minimal graphs $\{M\}$. The minimization procedure is completed if no sets of new vertices at the expansion stage lead to a successful iteration.

Next, we will take a closer look at each stage of the algorithm. The calculations at each step of both phases of the reduction are organized as follows: first, a list of vertex sets to be checked is compiled, then the entire list is checked. Thus, on the one hand, the task of obtaining such a list arises, but on the other hand, the possibility of parallel processing of the list items is immediately available.

\subsection{Expansion}

The expansion of accumulation graph $A$ is performed by filling in the partially occupied base orbits of graph $A$ and (or) adding new vertices from reserve graph $R$.

Graph $R$ is formed from base graph $B$ and is updated after several successful iterations of the minimization procedure. A base orbit is included in the graph $R$ if, when added to the graph $A$, the degree of some of its vertices is at least 4. (Note that the added base orbit generally splits into several orbits of graph $A$.) New orbits are added in order of decreasing vertex degree.

Our experience shows that during expansion, priority should be given to partially filled base orbits, and especially to orbits with a small number of vertices (6 or 12). The number of added vertices is a compromise between the calculation time (or the ability to calculate something at all) and the convergence of the search result to a global minimum.

We used the base graph $B=\oplus^4H^2$, choosing orbits that fall completely or partially on a disk of radius $r = 2$ (expanding the base graph and increasing the radius of the disk did not lead to new results in our case). The corresponding base orbits $(a, b, c, d)$ satisfy the relations
$$(a-b\sqrt{33})^2+(c\sqrt{3}-d\sqrt{11})^2\le(12 r)^2, \;\;a,b,c,d \in |2\mathbb{Z}|, \;a-b+c+d \in 4\mathbb{Z}.$$
There are only about 400 such orbits, of which about 40 are included in R. We noticed that all useful orbits satisfy the relation $abcd=0$.

Here we can consider the concept of an optimal minimization \textit{strategy}, i.e. a sequence of operations providing the fastest convergence to a global minimum. The universal optimal strategy is unknown. Usually a good strategy becomes clear only after a minimal graph is obtained.

It is wise to try different strategies and initial graphs. We used the following strategies, which give good results:
\begin{itemize}
\itemsep=0pt
    \item 
sequential addition of new base orbits in small portions
    \item 
filling and (or) adding a relatively large number of orbits, but with the subsequent removal of new vertices with a small final degree
    \item 
preliminary cycle of quick verification (to a small hyperedge degree) of the entire list of new vertices (with the addition of one or two orbits at each iteration) and allocation for further study of the most promising orbits by numerical indicators (see 5.10)
\end{itemize}

The concept of strategy also includes the choice of the general form of the graph based on guesses about its optimal symmetry or type. For example, one can expect that for a triples a good approximation will be provided by graphs with symmetries of order 3 or 6; for pairs with symmetries of order 2 or 4.

\subsection{Reduction}
The reduction stage includes the compilation of an ordered list of all (non-isomorphic) subsets of vertices of the graph $W$ containing at least one vertex of each hyperedge of $Y$. Usually it is more convenient to operate with sets of deleted vertices, instead of preserved ones, since they are much smaller. The procedure for constructing such a list may be as follows.

Using hyperedges of degree 2, we find all its maximal independent sets. To each such set we add all vertices (called free) that are not included in the hyperedges of degree 1 and 2. Next, we check each formed set for the presence of hyperedges of degree 3 and above. If the set contains some hyperedge of degree $n$, then such a set is divided into $n$ subsets, in each of which one of the vertices of this hyperedge is removed. The found sets are ordered by decreasing number of vertices. All sets that, when removed from $W$, give a graph with a number of vertices that exceeds the order of the current minimal graph $M$, are discarded.

Here, an \textit{independent vertex set} is a set of vertices such that no two vertices in the set are connected by an edge. A \textit{maximal} independent vertex set is an independent vertex set that is not a subset of any larger independent vertex set. It should be distinguished from a \textit{maximum} independent vertex set, which is an independent vertex set containing the largest possible number of vertices for a given graph.

\subsection{Accelerators}
To speed up the computing, a set of hardware and software \textit{accelerators} is used. Accelerators are an essential part of our method. So in our case, the total speed-up from using various accelerators reached about 100-fold (excluding SAT accelerators). The following accelerators gave the greatest effect:
\begin{itemize}
\itemsep=0pt
    \item 
parallel computing
    \item 
use of symmetries of input and output graphs
    \item 
one-time execution of repeated operations
    \item 
use of dependence of SAT check time on the result
    \item 
SAT accelerators (symmetry breaking and vertex ordering)
\end{itemize}
We now give a brief description of these accelerators.

\paragraph{Parallel computing.}
The minimization procedure requires a large number of independent similar checks carried out on precomputed lists of graphs. Therefore, it allows very efficient parallelization both between the kernels of the same processor and between different computers in the cluster. The degree of acceleration approaches the total number of used processing cores.

\paragraph{Using symmetries.}
The acceleration effect can be obtained by increasing the symmetry group order of both the input graph $W$ (due to thinning the list of checked graphs) and the output graph $M$ (due to the use of orbits instead of individual vertices).

The input and output symmetry can be changed by adding \textit{hanging vertices} to the graph $W$, that is vertices of degree 1, which can be taken as part of the input and output companion graphs $C$. The hanging vertices are those connected to a high-degree vertex of graph $W$, that are not deleted during reduction. The input and output symmetry groups are selected so as to form a quotient group. The sets of vertex permutations for the input and output symmetries of the graph $W$ are determined. Based on them the set of permutations of the quotient group is formed, which determines all automorphisms of $W$ in terms of output orbits.

Subsequently, from the list of subgraphs of the graph $W$ that require verification at the reduction stage, all isomorphic copies are deleted. This can be done by first replacing each subgraph in the list with the first element of the sorted list of all its isomorphic copies obtained by means of known set of permutations, and then removing the duplicated subgraphs from the list.

In addition to speeding up the computation, the use of symmetries allows us to expand the applicability of the method to larger graphs. At the initial iterations of the minimization procedure, one can operate with more symmetric graphs, then gradually lower the order of the input and output symmetries. The following undesirable “overcompression” effect was noticed on this path: the minimal asymmetric graph is usually not a subgraph of the minimal symmetric graph (this is especially noticeable for $S$-subgraphs). The problem of escaping the local minimum can be addressed by increasing the number of new vertices introduced in the expansion stage.

\paragraph{One-time execution of repeated operations.}
Since at the reduction stage all the graphs being checked are subgraphs of the same graph $W$, it is possible to reduce the preparation time of the SAT task by extracting the common part of the CNF-formula (see \ref{sat}). This part is calculated once, and immediately before the check, a small variable part corresponding to a certain subgraph is added to it. In our implementation, the clauses corresponding to all edges of the graph $W\cup C$ are placed in the common part, and those corresponding to surviving vertices are placed in the variable part.

\paragraph{Accelerated hypergraph construction.}
The time taken by the SAT solver depends significantly on the result of the test. So on average to get the result {\tt True} (meaning that 4 colors are enough for coloring) it takes an order of magnitude less time than to get the result {\tt False} (4 colors are not enough). In our case, for M-type graphs, the SAT solver need about 1 and 10 seconds respectively. This fact can be used to speed up calculations in the first phase of reduction, in which {\tt False} is a much more frequent test result (in the second phase, on the contrary, {\tt True} occurs much more frequently).

Various implementations of this idea are possible. In practice, the following simple “8-4-2-1” algorithm gives good results and on average leads to a 2-fold acceleration, despite an increase in the total number of checks. The sets of deleted vertices are first divided into groups of 8, and all vertices belonging to at least one set in the group are simultaneously deleted from the working graph. If the result is {\tt False}, then this means that each of the sets can be deleted individually. If the result is {\tt True} (which takes relatively little time), then the sets are grouped in fours, the test is repeated, then in pairs if needed, and finally individually.

\subsection{SAT encoding}
\label{sat}
A task is supplied to a SAT solver as a so-called \textit{propositional formula} in \textit{conjunctive normal form}, CNF. A formula includes the following components: variables, literals, clauses. A \textit{variable} describes some elementary statement, a \textit{literal} is some variable $v$ or its negation $\overline{v}$, a \textit{clause} is a union of literals by logical {\tt or} ($\vee$), and a \textit{formula} is a union of clauses by logical {\tt and} ($\wedge$). A formula is \textit{satisfiable} if it is possible to set all used variables to {\tt True} or {\tt False} in such a way that the overall formula simplifies to {\tt True}, i.e. each clause contains at least one literal that is {\tt True}.

To check the 4-colorability of an $n$-vertex graph $G=(V,E)$, $4n$ variables $v_{i, k}$ are introduced, meaning that the vertex $v_i$ has color $k$. A set of clauses includes vertex clauses $ v_{i,1} \vee v_{i,2} \vee \ldots \vee v_{i,k}:\forall v_i \in V$ (one per vertex, meaning that the vertex must be colored), and edge clauses $\overline{v_{i,k}} \vee \overline{v_{j,k}}: \forall \left\{v_i,v_j\right\} \in E$ (four per edge, meaning that adjacent vertices must have different colors).

To speed up solving, so-called \textit{symmetry breaking} clauses $v_{i,i}: \forall v_i \in Q$ are added to the formula, which assign different colors to the vertices of a specific 3-clique $Q\subset V$.
The choice of $Q$ significantly affects the speed of solving. For M-type graphs, good results are obtained by using the vertices of $L$ with coordinates $(0, 0)$, $(1, 0)$, $(1/2, \sqrt3/2)$. The same clique is used for centered non-mono-triples; moreover, a further 10-fold or so acceleration is possible for a triple with side $\sqrt3$, by means of combining symmetry-breaking with a companion graph $C$, and assigning colors to five vertices at once (vertices of non-mono-triple $(1, 0)$, $(-1/2, \sqrt3/2)$, $(-1/2, -\sqrt3/2)$ have identical color). For a triple with side $1/\sqrt3$ such a trick does not work, and an attempt to assign colors to four vertices slows down solving. For a mono-pair with length $8/3$ and coordinates $(4/3,0)$, $(-4/3,0)$, good results arise from choosing the vertices with coordinates $(4/3, 0)$, $(1/3, 0)$, $(5/6, \sqrt3/2)$.

Additionally, to block non-mono pairs, triples, and so on by a corresponding $m$-vertex mono-set one can use a cyclic chain of implications $v_{1,k} \Rightarrow v_{2,k} \Rightarrow ...\Rightarrow v_{m,k} \Rightarrow v_{1,k}$, meaning "if some vertex has given color, then the next vertex in the chain must have this color too", and so for each of $k$ colors. Implication $v_{i,k} \Rightarrow v_{j,k}$ corresponds to the clause $\overline{v_{i,k}} \vee v_{j,k}$. This approach requires fewer clauses than using graph $C$.

The order of vertices and, accordingly, the order of clauses, apparently also affects the speed of solving. We order the vertices by increasing the Euclidean distance from graph center, which gives good results.

\subsection{Rough reduction}
\label{rough}
We refer to the above method as the \textit{fine search}. Before using it, one should use a \textit{rough reduction} method to obtain an acceptable initial graph. Rough reduction methods include:
\begin{itemize}
\itemsep=0pt
    \item 
choice of graph type and base subgraphs
    \item 
clipping vertices beyond a certain Euclidean radius (called \textit{trimming})
    \item 
iterative removal of small-degree vertices
    \item
fixing vertices whose removal leads to loss of the key property
    \item 
reduction according to the orbits of the base graph (first of all, the orbits with the largest number of vertices are deleted)
    \item 
selection of graphs by some numerical indicators (the number of free vertices, the number of hyperedges of some small degrees, the order of the maximum independent set for hyperedges of degree 2)
\end{itemize}

In the latter case, the reduction is not carried out completely; instead, numerical characteristics are determined for a number of graphs by which the most promising graphs are selected for the subsequent reduction. This approach was used for S-subgraphs: first, from various subgraphs of the 187-vertex graph with the smallest number of hyperedges of degree up to 3, several 169-vertex graphs with symmetry of order 3 were selected, after which reduction was performed with output symmetry of order 1 and hyperedges up to degree 4. The total number of checks was slightly more than $10^{5}$ (compared to $10^{35}$ for brute force), which took slightly more than a day on a standard laptop.

\subsection{Limitations}
The maximum number of redundant vertices or orbits depends on the graph and the capabilities of the computer, and is usually limited to several tens. This is enough to quite effectively minimize 5-chromatic graphs. Taking into account the use of symmetries, the number of vertices is limited to about 1000.

\section{Results}
\label{results}

\renewcommand{\bottomfraction}{0.7}
\renewcommand{\textfraction}{0.0}

\paragraph{Summary.} The main results are summarized in table \ref{tres}. The columns indicate the type of graph, output symmetry order of a set of minimal graphs $\{M\}$ (and minimal graph $M$ with smallest edge number in parentheses) taking into account the connection edges, the number of vertices of minimal ($M$) and accumulative ($A$) graphs, the number of edges of graphs in $\{M\}$, and the number of non-isomorphic graphs in $\{M\}$. 

\begin{table}[!b]
{
\caption{Main results}
\label{tres}
\smallskip

{
\centering
\small
\begin{tabular}{c|c|c|c|c|c|c}
\hline
\multicolumn{2}{c|}{graph type} &\!sym\! &$M$ &edges &$A$ &set	\\
\hline
\hline
M-type graph	&M6A	&1	    &509	&2442	&		\\
	            &	    &1	    &510	&2502	&		\\
\hline
$L$-subgraph    &M6A	&1\,(2) &374	&1860-1884	&412  &1092	\\
	            &   	&2	    &374	&1868-1884	&406  &78	\\
	            &   	&1\,(2)	&375	&1920-1930	&403  &282	\\
            	&   	&3  	&376	&1890-1947	&451  &96	\\
\cline{2-7}
            	&M6B	&2	    &374	&1864-1880	&406  &78	\\
            	&   	&1\,(2)	&375	&1916-1926	&403  &66	\\
\cline{2-7}
            	&M6C	&2	    &374	&1872		&400  &76	\\
\hline
$S$-subgraph    &M6A    &1	    &136	&564-567	&166  &8	\\
\cline{2-7}
            	&M6B	&1	    &141	&594-595	&172  &4	\\
\cline{2-7}
            	&M6C	&1	    &150	&639-651	&167  &28	\\
\hline
mono-pair   &$8/3$	    &1\,(4)	&367	&1822		&375  &27	\\
\cline{2-7}
            &$8/\sqrt3$	&4	    &421	&2094-2106	&433  &7	\\
\hline
non-mono-pair 	&$3$    &4      &214	&977    	&226  &2	\\
\cline{2-7}
             	&$7/3$  &2\,(4) &319	&1524-1536	&343  &34	\\
\cline{2-7}
             	&$5/3$  &2      &315	&1498-1504	&335  &17	\\
\cline{2-7}
             	&$1/3$  &4      &312	&1502   	&312  &1	\\
\cline{2-7}
	   &$\sqrt{11/3}$	&2      &308	&1522-1524	&310  &2	\\
\hline
non-mono-triple &$\sqrt7$  &1   &159	&646-657	&180  &51	\\
\cline{2-7}
                &$\sqrt5$  &3   &265	&1293   	&265  &1	\\
\cline{2-7}
                &$\sqrt3$  &1   &221	&1006-1013	&238  &120  \\
\cline{2-7}
            &$\sqrt{5/3}$  &3   &250	&1116   	&250  &1	\\
\cline{2-7}
            &$\sqrt{1/3}$  &1   &265	&1246-1258	&289  &1020	\\
\hline
\end{tabular}

}
}
\end{table}


\begin{figure}[tb]
\centering
\includegraphics[scale=0.33]{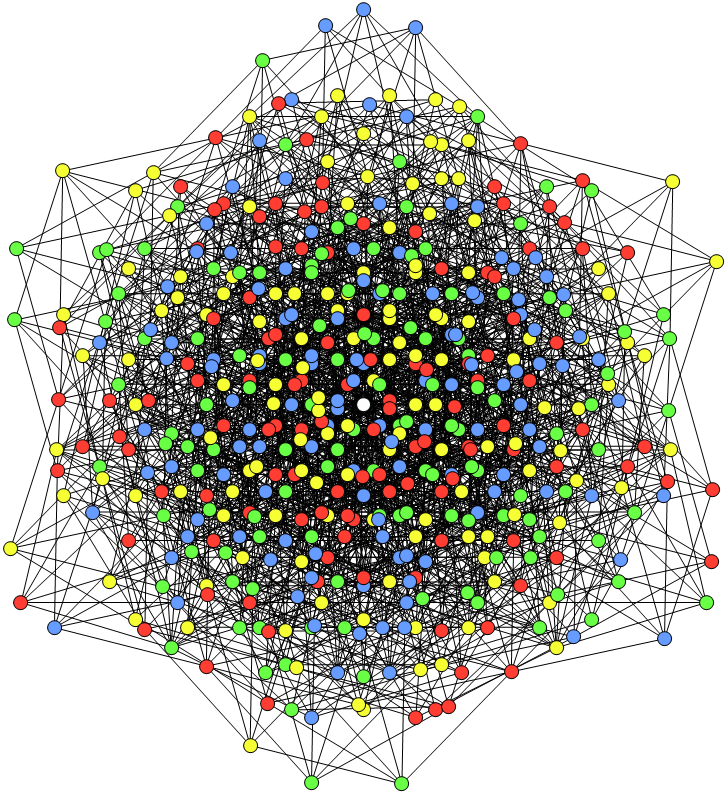}
\caption{5-chromatic graph with 509 vertices and 2442 edges.}
\label{g509}
\end{figure}

\begin{figure}[!p]
\centering
\includegraphics[scale=0.33]{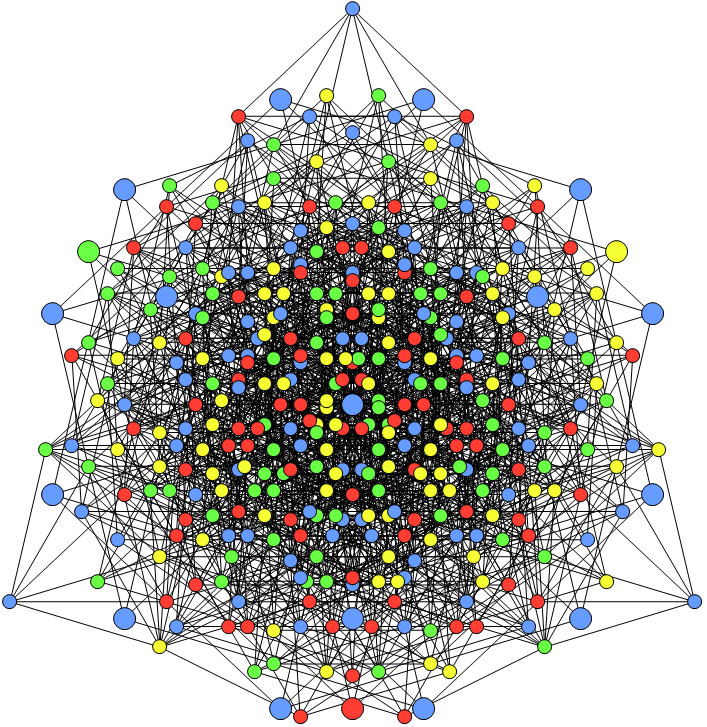}
\caption{Large subgraph $L_{374}$ with 374 vertices and 1860 edges.}
\label{l374}
\end{figure}

\begin{figure}[!p]
\centering
\includegraphics[scale=0.33]{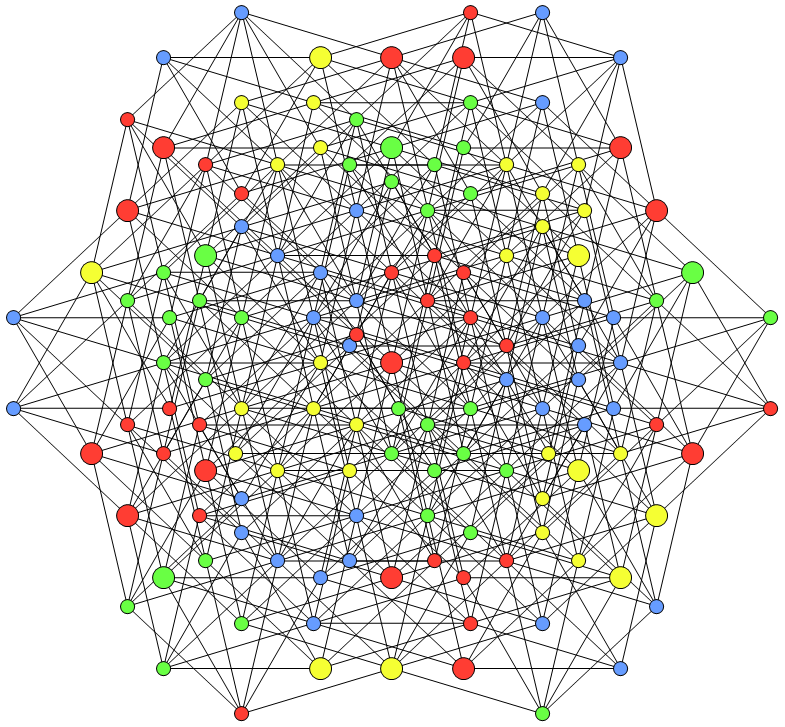}
\caption{Small subgraph $S_{136}$ with 136 vertices and 564 edges. Vertices, that can be used for connection to a large subgraph, are enlarged.}
\label{s136}
\end{figure}

The main attention was paid to graphs of type M. The record-holder is a graph with 509 vertices and 2442 edges (Fig.\ref{g509}), comprising a large subgraph with 374 vertices and 1860 edges (Fig.\ref{l374}) and a small subgraph with 136 vertices and 564 edges (Fig.\ref{s136}). For comparison we show the previous 510-vertex record obtained in the competition with Heule, and 375- and 376-vertex variants obtained for different orders of output symmetry. 
In addition, we discovered small graphs for a mono-pairs with a distance $8/3$, and $8/\sqrt3$, non-mono-pairs with distance $3$, $7/3$, $5/3$, $1/3$, and $\sqrt{11/3}$, and non-mono-triples with a side $\sqrt7$, $\sqrt5$, $\sqrt3$, $\sqrt{5/3}$, and $1/\sqrt3$.
Lists of the vertices of some graphs can be found on the Polymath project website \cite{par}.

\paragraph{M-type graphs.}
Surprisingly, comparable results were obtained for all three subtypes M6 of $L$-subgraphs. The subtypes M6A and M6B do not differ at all (except for the auxiliary orbit $(0, 0, 2, 6)$); the subtype M6C has a slightly different set of minimal graphs $\{M\} $. Moreover, for all subtypes the same composition of base orbits can be given to obtain a graph with the least number of both vertices and edges (the number of edges itself is different).

A completely different picture was observed for $S$-subgraphs, where subtype M6A is easily the best, and subtype M6C gives very weak results. (In fairness, we must note that subtypes M6B and M6C received much less attention, and it can be expected that the results will slightly converge.)

To minimize the $L$-subgraph we used an iterative procedure with alternating stages of reduction and adding a small number of orbits, and an output symmetry of order 3. If we start with the $L$-subgraph of Heule’s 529-vertex graph \cite{heu529}, then the path to further minimization opens after adding small-order orbits $(0, 0, 12, 0)$, $(0, 0, 2, 2)$, as well as the orbits $(4, 0, 10, 2)$, $(10, 0, 2, 4)$. Thereby we obtained a set of minimal 376-vertex graphs. Then, the corresponding 451-vertex accumulative graph was reduced with output symmetry of order 1, and a set of 375-vertex graphs was obtained. Note that none of the 376-vertex graphs is an overgraph of any of the 375-vertex graphs, and the union of the 375-vertex graphs is a 403-vertex graph.
The final results were obtained after changing the output symmetry from 3 to 2,  adding orbits $(0, 0, 0, 4)$ and $(0, 0, 6, 6)$, and filling orbit $(4, 0, 0, 0)$.

Reduction of $S$-subgraphs is briefly discussed in \ref{rough}. Here we hardly improved the result of Heule \cite{heu510} at all; only two edges were removed.

\paragraph{Orbit-filling tables.}
Below are the images and \textit{orbit-filling} tables of some graphs that are smallest at the moment. The tables list the base orbits set, the maximum vertex degree for each base orbit, and the number of vertices for minimal graph $M$, for all graphs of the set of minimal graphs $\{M\}$ and for the accumulative graph $A=\bigcup\{M\}$. 
The base orbits set is formed from those orbits of the base graph (with a symmetry group of order 24) that include at least one vertex from $A$.

Table \ref{tcmp} shows the differences in orbit-filling for $L$-subgraphs of some minimal graphs shown in Table \ref{tres}. For comparison, we include the $L$-subgraphs with 375, 382, and 393 vertices obtained by Heule (corresponding to his graphs with 510, 517, and 529 vertices).
	

\begin{table}[!h]
{
\caption{Comparison of orbits filling of some $L$-subgraphs. 
}
\label{tcmp}

{
\centering
\small
\begin{tabular}{c|c|c|c|c|c|c|c|c}
\hline
            & \multicolumn{4}{c|}{$M$}  & \multicolumn{4}{c}{$\{M\}$}  \\
\cline{2-9}
orbit      &374 &375 &382 &393 & 376   & 375  & 374  & 374  \\
\cline{2-9}
 & \multicolumn{2}{c|}{M6A\&B}  & \multicolumn{2}{c|}{M6A}  & 
 M6A  & \multicolumn{2}{c|}{M6A\&B}  & M6C  \\
\hline
\hline
\gr (0,0,0,0)  &\gr 1  &\gr 1  &\gr 1  &\gr 1  &\gr 1  &\gr 1  &\gr 1  &\gr 1  \\
\hline
 (4,0,0,0)  & 6  & 4  & 6  & 6  & 3     & 4    & 6    & 6    \\
\gr (0,0,4,0)  &\gr 6  &\gr 6  &\gr 6  &\gr 6  &\gr 6  &\gr 6  &\gr 6  &\gr 6  \\
\gr (12,0,0,0) &\gr 6  &\gr 6  &\gr 6  &\gr 6  &\gr 6  &\gr 6  &\gr 6  &\gr 6  \\
 (0,0,0,4)  & 3  & 0  & 0  & 0  & 0     & 0    & 3    & 3    \\
\gr (0,0,8,0)  &\gr 6  &\gr 6  &\gr 6  &\gr 6  &\gr 6  &\gr 6  &\gr 6  &\gr 6  \\
 (0,0,12,0) & 6  & 6  & 0  & 0  & 6     & 6    & 6    & 6    \\
\hline
 (0,0,2,2)  & 8  & 6  & 0  & 0  & 6     & 6    & 6-8  & 8    \\
\gr (6,2,0,0)  &\gr 12 &\gr 12 &\gr 12 &\gr 12 &\gr 12 &\gr 12 &\gr 12 &\gr 12 \\
\gr (2,0,0,2)  &\gr 12 &\gr 12 &\gr 12 &\gr 12 &\gr 12 &\gr 12 &\gr 12 &\gr 12 \\
\gr (6,0,0,2)  &\gr 12 &\gr 12 &\gr 12 &\gr 12 &\gr 12 &\gr 12 &\gr 12 &\gr 12 \\
 (0,0,6,6)  & 6  & 0  & 0  & 0  & 0-3   & 0    & 6    & 6    \\
\gr (10,0,0,2) &    6  &    6  &    6  &    6  &    6  &    6  &    6  &    6  \\
\gr (0,2,2,0)  &\gr 12 &\gr 12 &\gr 12 &\gr 12 &\gr 12 &\gr 12 &\gr 12 &\gr 12 \\
\gr (4,0,0,4)  &\gr 12 &\gr 12 &\gr 12 &\gr 12 &\gr 12 &\gr 12 &\gr 12 &\gr 12 \\
\gr (14,0,0,2) &    6  &    6  &    6  &    6  &    6  &    6  &    6  &    6  \\
\gr (0,2,6,0)  &\gr 12 &\gr 12 &\gr 12 &\gr 12 &\gr 12 &\gr 12 &\gr 12 &\gr 12 \\
 (8,0,0,4)  & 0  & 0  & 6  & 6  & 0     & 0    & 0    & 0    \\
 (0,0,2,6)  & 6  & 6  & 6  & 6  & 6     & 6    & 6    & 6    \\
\gr (6,0,10,0) &\gr 12 &\gr 12 &\gr 12 &\gr 12 &\gr 12 &\gr 12 &\gr 12 &\gr 12 \\
\gr (6,0,0,6)  &\gr 12 &\gr 12 &\gr 12 &\gr 12 &\gr 12 &\gr 12 &\gr 12 &\gr 12 \\
\gr (0,4,4,0)  &\gr 12 &\gr 12 &\gr 12 &\gr 12 &\gr 12 &\gr 12 &\gr 12 &\gr 12 \\
\hline
\gr (2,0,4,2)  &\gr 24 &\gr 24 &\gr 24 &\gr 24 &\gr 24 &\gr 24 &\gr 24 &\gr 24 \\
 (12,2,2,0) & 12 & 12 & 12 & 15 & 12    & 12   & 12   & 12   \\
 (2,0,6,4)  & 2  & 4  & 6  & 0  & 0-12  & 0-4  & 0-8  & 0-8  \\
\gr (4,0,2,2)  &\gr 24 &\gr 24 &\gr 24 &\gr 24 &\gr 24 &\gr 24 &\gr 24 &\gr 24 \\
 (4,0,6,2)  & 12 & 12 & 12 & 15 & 12    & 12   & 12   & 12   \\
\gr (2,0,8,2)  &    12 &    12 &    12 &    12 &    12 &    12 &    12 &    12 \\
 (4,0,4,4)  & 0  & 0  & 3  & 6  & 0-3   & 0    & 0    & 0    \\
\gr (8,0,2,2)  &\gr 24 &\gr 24 &\gr 24 &\gr 24 &\gr 24 &\gr 24 &\gr 24 &\gr 24 \\
 (8,0,6,2)  & 12 & 12 & 12 & 6  & 12    & 12   & 12   & 12   \\
 (6,2,4,0)  & 12 & 24 & 24 & 21 & 18-24 & 24   & 12   & 12   \\
 (2,0,2,4)  & 8  & 10 & 9  & 17 & 9-12  & 10   & 8-10 & 8    \\
\gr (10,0,4,2) &    12 &    12 &    12 &    12 &    12 &    12 &    12 &    12 \\
 (8,0,4,4)  & 12 & 12 & 12 & 12 & 0-12  & 12   & 12   & 12   \\
 (12,2,6,0) & 0  & 0  & 0  & 6  & 0-12  & 0    & 0    & 0    \\
 (4,0,10,2) & 0  & 6  & 0  & 0  & 0-12  & 2-6  & 0-4  & 0    \\
 (10,0,2,4) & 10 & 8  & 0  & 0  & 0-12  & 8-12 & 4-12 & 4-12 \\
 (6,0,4,6)  & 0  & 0  & 3  & 6  & 0     & 0    & 0    & 0    \\
 (6,2,8,0)  & 12 & 6  & 24 & 24 & 0-12  & 6-10 & 8-12 & 12   \\
\gr (14,0,2,4) &    12 &    12 &    12 &    12 &    12 &    12 &    12 &    12 \\
\hline
\end{tabular}

}
}
\end{table}

Base orbits can be divided into \textit{stable} (with the same, nonempty, vertex set for all minimal graphs) and \textit{unstable} (with some minimal graphs having more vertices than others). Stable orbits can be divided into \textit{filled} (with the maximum available number of vertices for a given orbit) and \textit{unfilled}. The latter usually include exactly half or a quarter of the maximum number of vertices.

The orbits are grouped by the number of vertices (which for a filled orbit is 1, 6, 12 or 24) and listed in increasing order of the smaller of the two possible Euclidean radii. Filled orbits  are highlighted in gray. Stable orbits are highlighted in gray in first column of Tables 2 and 3.

The orbit-filling tables are usually more informative than the graph image, and allow one to compare different graphs and to devise directions of further attack in the minimization problem.

\paragraph{Hardware, software and execution time. }
All calculations were performed in {\tt Mathemetica 10} on a laptop with processor Intel Core i7-2670QM, 2.2GHz, 4 cores/8 threads, 6 GB RAM.
We used hyperedges of degree 2 or 3 for $L$-subgraph, and of degrees 3 or 4 for $S$-subgraph and non-mono-triples. Some pairs and triples (with a final number of vertices less than 250) made it possible to effectively use hyperedges of degrees up to 5.
The average time for a single check with output result {\tt False} (unsatisfiable) was about 10 s for M-type graphs, and about 0.2~s for non-mono-triples with side $\sqrt3$. The typical time to minimize $L$- or $S$-subgraphs was several tens of hours.


\newpage

\begin{table}[!p]
{
\caption{Filling orbits of $L$-subgraph M6A or M6B, $M=G_{374}$, $A=G_{412}$}
\label{tl}
\smallskip
\small
\begin{tabular}{cc}
\begin{tabular}{c|c|c|c|c}
\hline
 orbit      & deg & $M$  & $\{M\}$ & $A$  \\
\hline
\hline
\gr (0,0,0,0)  & 30 &\gr  1 &\gr 1  &\gr 1  \\
\hline
    (4,0,0,0)  & 8  &\gr  6 &\gr 6  &\gr 6  \\
\gr (0,0,4,0)  & 26 &\gr  6 &\gr 6  &\gr 6  \\
\gr (12,0,0,0) & 23 &\gr  6 &\gr 6  &\gr 6  \\
    (0,0,0,4)  & 6  &     3 &    3  &    3  \\
\gr (0,0,8,0)  & 14 &\gr  6 &\gr 6  &\gr 6  \\
    (0,0,12,0) & 6  &\gr  6 &\gr 6  &\gr 6  \\
\hline
    (0,0,2,2)  & 12 &     8 &  6-8  &    8  \\
\gr (6,2,0,0)  & 17 &\gr 12 &\gr 12 &\gr 12 \\
\gr (2,0,0,2)  & 19 &\gr 12 &\gr 12 &\gr 12 \\
\gr (6,0,0,2)  & 9  &\gr 12 &\gr 12 &\gr 12 \\
    (0,0,6,6)  & 6  &     6 &    6  &    6  \\
\gr (10,0,0,2) & 15 &     6 &    6  &    6  \\
\gr (0,2,2,0)  & 25 &\gr 12 &\gr 12 &\gr 12 \\
\gr (4,0,0,4)  & 11 &\gr 12 &\gr 12 &\gr 12 \\
\gr (14,0,0,2) & 9  &     6 &    6  &    6  \\
\gr (0,2,6,0)  & 13 &\gr 12 &\gr 12 &\gr 12 \\
    (0,0,2,6)  & 6  &     6 &    6  &    6  \\
\gr (6,0,10,0) & 11 &\gr 12 &\gr 12 &\gr 12 \\
\hline
\end{tabular}
&
\begin{tabular}{c|c|c|c|c}
\hline
 orbit      & deg & $M$  & $\{M\}$ & $A$  \\
\hline
\hline
\gr (6,0,0,6)  & 10 &\gr 12 &\gr 12 &\gr 12 \\
\gr (0,4,4,0)  & 4  &\gr 12 &\gr 12 &\gr 12 \\
\hline
\gr (2,0,4,2)  & 17 &\gr 24 &\gr 24 &\gr 24 \\
    (12,2,2,0) & 14 &    12 &    12 &    12 \\
    (2,0,6,4)  & 7  &     2 &  0-8  &    12 \\
\gr (4,0,2,2)  & 16 &\gr 24 &\gr 24 &\gr 24 \\
    (4,0,6,2)  & 11 &    12 &    12 &    12 \\
\gr (2,0,8,2)  & 10 &    12 &    12 &    12 \\
\gr (8,0,2,2)  & 17 &\gr 24 &\gr 24 &\gr 24 \\
    (8,0,6,2)  & 13 &    12 &    12 &    12 \\
    (6,2,4,0)  & 13 &    12 &    12 &    12 \\
    (2,0,2,4)  & 10 &     8 &  8-10 &    12 \\
\gr (10,0,4,2) & 12 &    12 &    12 &    12 \\
    (8,0,4,4)  & 8  &    12 &    12 &    12 \\
    (4,0,10,2) & 6  &     0 &   0-4 &    12 \\
    (10,0,2,4) & 7  &    10 &  4-12 &    12 \\
    (6,2,8,0)  & 10 &    12 &  8-12 &    12 \\
\gr (14,0,2,4) & 8  &    12 &    12 &    12 \\
 &&&&\\
\hline
\end{tabular}
\end{tabular}
}
\end{table}

\newpage

\begin{table}[!p]
{
\caption{Filling orbits of minimal $S$-subgraphs of subtype M6. Corresponding graphs are of order 136 (M6A), 141 (M6B), and 150 (M6C).}
\label{ts}
\smallskip

{
\centering
\small
\begin{tabular}{c|c|c|c|c|c|c|c}

\hline
  &  & \multicolumn{3}{c|}{$M$}  & \multicolumn{3}{c}{$\{M\}$}  \\
\cline{3-8}
 orbit      &deg & M6A & M6B & M6C & M6A & M6B & M6C  \\
\hline
\hline
 (0,0,0,0)  & 12 &\gr 1  &\gr 1  &\gr 1  &\gr 1  &\gr 1  &\gr 1  \\
\hline
 (0,0,4,0)  & 16 &\gr 6  &\gr 6  &\gr 6  &\gr 6  &\gr 6  &\gr 6  \\
 (0,0,8,0)  & 8  & 3     &    3  &    4  & 3     & 3     &    4  \\
 (0,4,0,0)  & 6  & 4     &    3  &    5  & 4-5   & 3     &   3-5 \\
\hline
 (6,2,0,0)  & 15 &\gr 12 &\gr 12 &\gr 12 &\gr 12 &\gr 12 &\gr 12 \\
 (0,2,2,0)  & 16 &\gr 12 &\gr 12 &\gr 12 &\gr 12 &\gr 12 &\gr 12 \\
 (0,2,6,0)  & 11 &\gr 12 &\gr 12 &\gr 12 &\gr 12 &\gr 12 &\gr 12 \\
 (0,0,2,6)  & 10 &\gr 12 &\gr 12 &\gr 12 &\gr 12 &\gr 12 &\gr 12 \\
 (6,0,0,6)  & 9  &\gr 12 &\gr 12 &\gr 12 &\gr 12 &\gr 12 &\gr 12 \\
 (0,4,4,0)  & 8  &\gr 12 &\gr 12 &\gr 12 &\gr 12 &\gr 12 &\gr 12 \\
\hline
 (12,2,2,0) & 6  & 3     &    4  &    8  & 3     & 4     &    8  \\
 (6,2,4,0)  & 11 & 14    &    15 &    13 & 14-16 & 15    & 12-16 \\
 (6,0,4,6)  & 8  & 9     &    13 &    17 & 7-9   & 13    & 16-19 \\
 (0,2,4,6)  & 7  &\gr 24 &\gr 24 &\gr 24 &\gr 24 &\gr 24 &\gr 24 \\
\hline
\end{tabular}

}
}
\end{table}

\newpage

\begin{table}[!p]
{
\caption{Filling orbits of $8/3$-mono-pair graph, $M=G_{367}$, $A=G_{375}$}
\label{t367}
\smallskip
\small
\begin{tabular}{cc}
\begin{tabular}{c|c|c|c|c}
\hline
 orbit      & deg & $M$  & $\{M\}$ & $A$  \\
\hline
\hline
 (0,0,0,0)  & 18 &\gr 1  &\gr 1  &\gr 1  \\
\hline
 (4,0,0,0)  & 28 &\gr 6  &\gr 6  &\gr 6  \\
 (8,0,0,0)  & 21 &\gr 6  &\gr 6  &\gr 6  \\
 (12,0,0,0) & 9  &\gr 6  &\gr 6  &\gr 6  \\
 (0,0,0,4)  & 14 &    4  &    4  &    6  \\
 (16,0,0,0) & 5  &    2  &    2  &    2  \\
\hline
 (10,2,0,0) & 11 &    2  &    2  &    4  \\
 (14,2,0,0) & 9  &    4  &    4  &    8  \\
 (0,0,2,2)  & 22 &\gr 12 &\gr 12 &\gr 12 \\
 (0,0,6,2)  & 16 &    10 &    10 &\gr 12 \\
 (0,0,4,4)  & 10 &    6  &    6  &    8  \\
 (2,0,0,2)  & 19 &\gr 12 &\gr 12 &\gr 12 \\
 (6,0,0,2)  & 19 &\gr 12 &\gr 12 &\gr 12 \\
 (2,2,0,0)  & 16 &\gr 12 &\gr 12 &\gr 12 \\
 (2,0,6,0)  & 17 &\gr 12 &\gr 12 &\gr 12 \\
 (10,0,0,2) & 13 &\gr 12 &\gr 12 &\gr 12 \\
 (4,0,8,0)  & 13 &    8  &    8  &    8  \\
\hline

\hline
\end{tabular}
&
\begin{tabular}{c|c|c|c|c}
\hline
 orbit      & deg & $M$  & $\{M\}$ & $A$  \\
\hline
\hline
 (2,0,4,2)  & 13 &\gr 24 &\gr 24 &\gr 24 \\
 (8,2,2,0)  & 10 &    8  &    8  &    8  \\
 (4,0,2,2)  & 16 &\gr 24 &\gr 24 &\gr 24 \\
 (4,0,6,2)  & 11 &\gr 24 &\gr 24 &\gr 24 \\
 (6,0,4,2)  & 15 &\gr 24 &\gr 24 &\gr 24 \\
 (6,0,6,4)  & 9  &    12 &    12 &    12 \\
 (4,2,2,0)  & 14 &    20 &    20 &    20 \\
 (8,0,2,2)  & 10 &    16 &    16 &    16 \\
 (8,0,6,2)  & 7  &    0  &   0-4 &    8  \\
 (6,0,8,2)  & 11 &    12 &    12 &    12 \\
 (2,0,2,4)  & 8  &    4  &    4  &    4  \\
 (10,0,4,2) & 8  &    12 &    12 &    12 \\
 (0,2,4,2)  & 8  &    12 &    12 &    12 \\
 (2,2,4,0)  & 11 &    12 &    12 &    12 \\
 (12,0,2,2) & 14 &\gr 24 &\gr 24 &\gr 24 \\
 (12,0,6,2) & 7  &    8  &   4-8 &    8  \\
 (4,2,6,0)  & 10 &    4  &    4  &    8  \\
\hline
\end{tabular}
\end{tabular}
}
\end{table}

\begin{figure}[!p]
\centering
\includegraphics[scale=0.33]{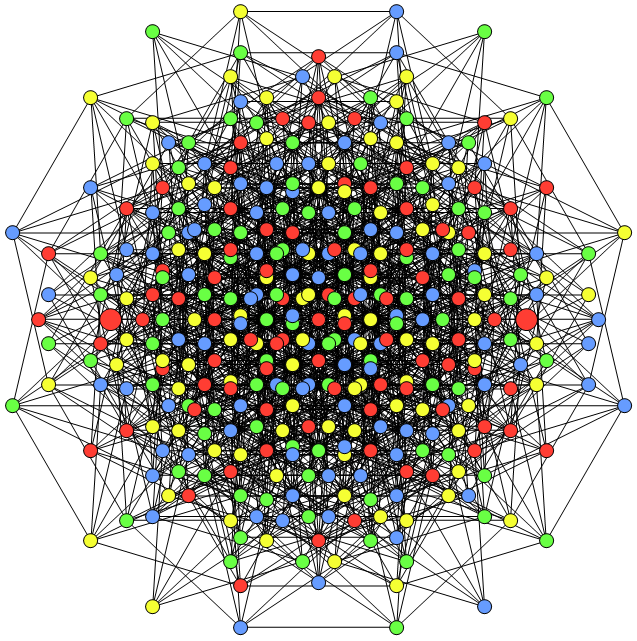}
\caption{The graph $G_{367}$ including monochromatic pair of length $8/3$ (corresponding vertices are enlarged), 367 vertices and 1822 edges}
\label{g367}
\end{figure}

\newpage

\begin{table}[!p]
{
\caption{Filling orbits of minimal graphs including non-mono-triple with side $\sqrt3$, $M=G_{221}$, $A=G_{238}$}
\label{t221}
\smallskip
\small
\begin{tabular}{cc}
\begin{tabular}{c|c|c|c|c}
\hline
 orbit      & deg & $M$  & $\{M\}$ & $A$  \\
\hline
\hline
 (0,0,0,0)  & 30 &\gr 1  &\gr 1  &\gr 1  \\
\hline
 (0,0,4,0)  & 24 &\gr 6  &\gr 6  &\gr 6  \\
 (12,0,0,0) & 17 &\gr 6  &\gr 6  &\gr 6  \\
 (0,0,8,0)  & 18 &\gr 6  &\gr 6  &\gr 6  \\
\hline
 (6,2,0,0)  & 15 &\gr 12 &\gr 12 &\gr 12 \\
 (2,0,0,2)  & 16 &    10 &    10 &\gr 12 \\
 (10,0,0,2) & 12 &    9  &    9  &    9  \\
 (0,2,2,0)  & 14 &\gr 12 &\gr 12 &\gr 12 \\
 (14,0,0,2) & 5  &    3  &    3  &    3  \\
 (0,2,6,0)  & 12 &\gr 12 &\gr 12 &\gr 12 \\
 (8,0,0,4)  & 7  &    3  &    3  &    3  \\
 (0,2,10,0) & 7  &    4  &   0-4 &    6  \\
\hline
\end{tabular}
&
\begin{tabular}{c|c|c|c|c}
\hline
 orbit      & deg & $M$  & $\{M\}$ & $A$  \\
\hline
\hline
 (2,0,4,2)  & 14 &    22 & 22-24 &\gr 24 \\
 (12,2,2,0) & 10 &    12 &    12 &    12 \\
 (4,0,2,2)  & 15 &\gr 24 &\gr 24 &\gr 24 \\
 (4,0,6,2)  & 11 &    14 & 14-18 &    18 \\
 (4,0,4,4)  & 7  &    6  &    6  &    6  \\
 (8,0,2,2)  & 11 &    14 & 12-14 &    18 \\
 (6,2,4,0)  & 14 &    21 &    21 &\gr 24 \\
 (2,0,2,4)  & 9  &    6  &    6  &    6  \\
 (10,0,4,2) & 7  &    6  &    6  &    6  \\
 (8,0,4,4)  & 8  &    6  &    6  &    6  \\
 (12,2,6,0) & 7  &    6  &    6  &    6  \\
 &&&&\\
\hline
\end{tabular}
\end{tabular}
}
\end{table}

\begin{figure}[!p]
\centering
\includegraphics[scale=0.33]{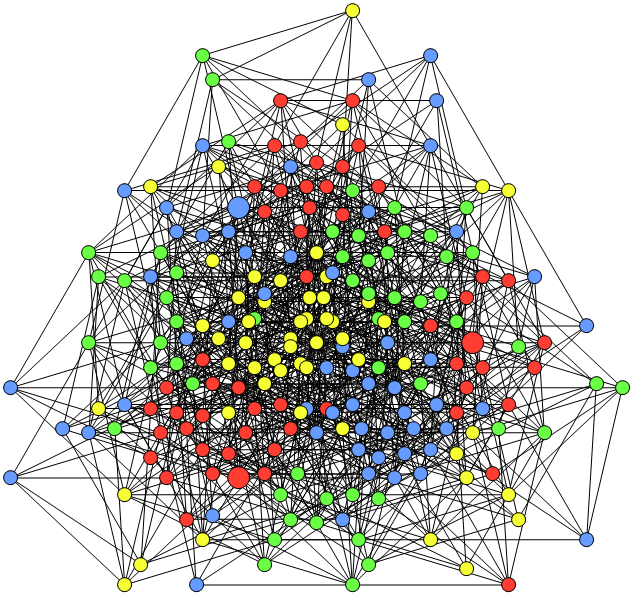}
\caption{The graph $G_{221}$ including non-monochromatic triple with side $\sqrt3$ (corresponding vertices are enlarged), 221 vertices and 1006 edges}
\label{g221}
\end{figure}

\newpage

\begin{table}[!p]
{
\caption{Filling orbits of minimal graphs including non-mono-triple with side $1/\sqrt3$, $M=G_{265}$, $A=G_{289}$}
\label{t265}
\smallskip
\small
\begin{tabular}{cc}
\begin{tabular}{c|c|c|c|c}
\hline
 orbit      & deg & $M$  & $\{M\}$ & $A$  \\
\hline
\hline
 (0,0,0,0)  & 30 &\gr 1  &\gr 1  &\gr 1  \\
\hline
 (4,0,0,0)  & 12 &    3  &    3  &    3  \\
 (0,0,4,0)  & 24 &\gr 6  &\gr 6  &\gr 6  \\
 (12,0,0,0) & 21 &\gr 6  &\gr 6  &\gr 6  \\
 (0,0,8,0)  & 12 &\gr 6  &\gr 6  &\gr 6  \\
\hline
 (0,0,2,2)  & 12 &    6  &    6  &    6  \\
 (6,2,0,0)  & 15 &\gr 12 &\gr 12 &\gr 12 \\
 (2,0,0,2)  & 20 &\gr 12 &\gr 12 &\gr 12 \\
 (6,0,0,2)  & 9  &    6  &    6  &    6  \\
 (10,0,0,2) & 12 &    6  &    6  &    6  \\
 (0,2,2,0)  & 18 &\gr 12 &\gr 12 &\gr 12 \\
 (4,0,0,4)  & 11 &    4  &    4  &    9  \\
 (14,0,0,2) & 11 &    6  &    6  &    6  \\
 (0,2,6,0)  & 10 &\gr 12 &\gr 12 &\gr 12 \\
 (8,0,0,4)  & 10 &    2  &    2  &    3  \\
\hline
\end{tabular}
&
\begin{tabular}{c|c|c|c|c}
\hline
 orbit      & deg & $M$  & $\{M\}$ & $A$  \\
\hline
\hline
 (2,0,4,2)  & 14 &    18 &    18 &    18 \\
 (12,2,2,0) & 12 &    12 &    12 &    12 \\
 (2,0,6,4)  & 8  &    0  &   0-2 &    6  \\
 (4,0,2,2)  & 17 &\gr 24 &\gr 24 &\gr 24 \\
 (4,0,6,2)  & 12 &    12 &    12 &    12 \\
 (2,0,8,2)  & 8  &    6  &    6  &    6  \\
 (4,0,4,4)  & 7  &    5  &   0-5 &    12 \\
 (8,0,2,2)  & 18 &\gr 24 &\gr 24 &\gr 24 \\
 (8,0,6,2)  & 12 &    12 &    12 &    12 \\
 (6,2,4,0)  & 13 &    12 &    12 &    12 \\
 (2,0,2,4)  & 10 &    12 &    12 &    12 \\
 (10,0,4,2) & 9  &    7  &  7-11 &    12 \\
 (12,0,2,2) & 5  &    3  &    3  &    6  \\
 (10,0,2,4) & 7  &    12 &    12 &    12 \\
 (6,2,8,0)  & 8  &    6  &    6  &    6  \\
\hline
\end{tabular}
\end{tabular}
}
\end{table}

\begin{figure}[!p]
\centering
\includegraphics[scale=0.33]{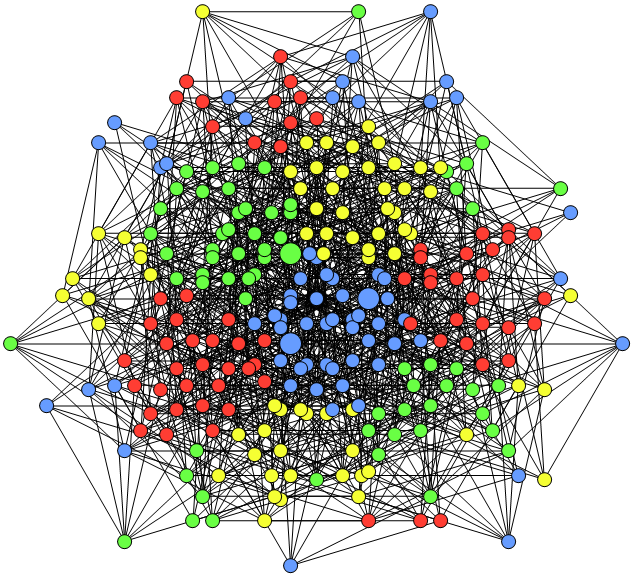}
\caption{The graph $G_{265}$ including non-monochromatic triple with side $1/\sqrt3$ (corresponding vertices are enlarged), 265 verices and 1246 edges}
\label{g265}
\end{figure}

\newpage

\newpage

\section{Conclusions}

We have introduced a new graph minimization method. This method can be useful in various minimization problems, where it is possible to effectively check whether a necessary property is preserved or not.

Initially, our method was conceived as not requiring special programs and large computing power, while providing a good local minimum. In this sense, our approach can be seen as an alternative to the approach of Heule \cite{heu553}. With regard to achieving record results, both approaches give comparable results after proper tuning. There are certain difficulties in comparing their computational efficiency. Heule estimated the total costs to compute his 529-vertex graph as 100\,000 CPU hours \cite{heu529}. Later he rated 1000 hours for a 510-vertex graph. We spent on the development of our approach roughly 1000 laptop hours, including the writing of programs and the study of different graphs. A 509-vertex graph can be found from scratch in about 100 laptop hours.

In contrast to the approach of Heule with randomization and extraction of the (first available) unsatisfiable core, our reduction method looks through all the solutions and finds a global minimum for a given graph (but cannot work with as many vertices). 
The success of reduction is largely determined by the choice of the initial graph. Our approach reduces the impact of this choice. This allowed us to find a more efficient set of orbits and move a little further.

As for the numerical results, they were obtained literally at the last minute. Not all options have been studied yet, which leaves the possibility of further progress.

It remains to be recalled that while we count the vertices, the enchanted princess is still wandering somewhere, with two extra heads, and is waiting for her hero.

\section{Acknowledgements}
I thank Aubrey de Grey for the opportunity to participate in the mind games, good support and text correction, Marijn Heule for preparing the sports ground, participants of the Polymath project for inexhaustible sea of knowledge that cannot be drunk in any way, and Alexander Soifer for his colorful book and generous offer, made in a Russian way, to insert me into the history.

\end{document}